\documentclass[reqno,12pt]{amsart}
\usepackage{amsmath, latexsym, amsfonts, amssymb, amsthm, amscd, color}
\usepackage{graphics,epsf,psfrag}
\DeclareMathOperator{\var}{\mathrm{Var}}

\numberwithin{equation}{section}

\newtheorem{theorem}{Theorem}[section]
\newtheorem{lemma}[theorem]{Lemma}
\newtheorem{proposition}[theorem]{Proposition}
\newtheorem{cor}[theorem]{Corollary}
\newtheorem{rem}[theorem]{Remark}
\newtheorem{definition}[theorem]{Definition}

\newcommand{\ind}{\mathbf{1}}

\newcommand{\R}{\mathbb{R}}
\newcommand{\Z}{\mathbb{Z}}
\newcommand{\N}{\mathbb{N}}
\renewcommand{\tilde}{\widetilde}
\renewcommand{\hat}{\widehat}

\DeclareMathSymbol{\leqslant}{\mathalpha}{AMSa}{"36} % nicer `smaller or equal'
\DeclareMathSymbol{\geqslant}{\mathalpha}{AMSa}{"3E} % nicer `larger or equal'
\DeclareMathSymbol{\eset}{\mathalpha}{AMSb}{"3F}     % nicer `emptyset'
\renewcommand{\leq}{\;\leqslant\;}                   % redef. of < or =
\renewcommand{\geq}{\;\geqslant\;}                   % redef. of > or =
\newcommand{\dd}{\,\text{\rm d}}             % a straight d for differentials
\newcommand{\bbL}{{\ensuremath{\mathbb L}} }

\newcommand{\bbR}{{\ensuremath{\mathbb R}} }

\newcommand{\ga}{\alpha}
\newcommand{\gb}{\beta}
\newcommand{\gga}{\gamma}            % \gg already exists...
\newcommand{\gd}{\delta}
\newcommand{\gep}{\varepsilon}       % \ge already exists...

\newcommand{\gG}{\Gamma}

\newcommand{\go}{\omega}

\newcommand{\gl}{\lambda}

\renewcommand{\ge}{\geq}
\renewcommand{\le}{\leq}
\makeatletter
\def\captionfont@{\footnotesize}
\def\captionheadfont@{\scshape}

\long\def\@makecaption#1#2{%
  \vspace{2mm}
  \setbox\@tempboxa\vbox{\color@setgroup
    \advance\hsize-6pc\noindent
    \captionfont@\captionheadfont@#1\@xp\@ifnotempty\@xp
        {\@cdr#2\@nil}{.\captionfont@\upshape\enspace#2}%
    \unskip\kern-6pc\par
    \global\setbox\@ne\lastbox\color@endgroup}%
  \ifhbox\@ne % the normal case
    \setbox\@ne\hbox{\unhbox\@ne\unskip\unskip\unpenalty\unkern}%
  \fi
  \ifdim\wd\@tempboxa=\z@ % this means caption will fit on one line
    \setbox\@ne\hbox to\columnwidth{\hss\kern-6pc\box\@ne\hss}%
  \else % tempboxa contained more than one line
    \setbox\@ne\vbox{\unvbox\@tempboxa\parskip\z@skip
        \noindent\unhbox\@ne\advance\hsize-6pc\par}%
\fi
  \ifnum\@tempcnta<64 % if the float IS a figure...
    \addvspace\abovecaptionskip
    \moveright 3pc\box\@ne
  \else % if the float IS NOT a figure...
    \moveright 3pc\box\@ne
    \nobreak
    \vskip\belowcaptionskip
  \fi
\relax
}
\makeatother
\def\writefig#1 #2 #3 {\rlap{\kern #1 truecm
\raise #2 truecm \hbox{#3}}}

\author{Hubert Lacoin}

\address{\noindent Universit\'e Paris 7, Math\'ematiques,
 case 7012, 2, Place Jussieu,
75251 Paris, France}

\email{lacoin@math.jussieu.fr}

\author{Gregorio Moreno}

\address{\noindent Universit\'e Paris 7, Math\'ematiques,
 case 7012, 2, Place Jussieu,
75251 Paris, France
{\rm and}
\noindent Facultad de Matem\'aticas, Pontificia Universidad Cat\'olica de Chile,
Vicu\~na Mackena 4860, Macul, Chile
}

\email{gregorio.random@gmail.com}

\title[Polymers on Hierarchical Lattices]{Directed Polymers on Hierarchical Lattices
with site disorder.}

\begin{document}

\maketitle

\begin{abstract}
We study a polymer model on hierarchical lattices very close to the one introduced and studied in
\cite{DGr, CD}. For this model, we prove the existence of free energy and derive the necessary and sufficient condition for which very strong disorder holds for all $\gb$, and give some accurate results on the behavior of the free energy at high-temperature. 
We obtain these results by using a combination of fractional moment method and change of measure over the environment to obtain an upper bound, and second moment method to get a lower bound. We also get lower bounds on the fluctuation exponent of $\log Z_n$, and study the infinite polymer measure in the weak disorder phase.
\end{abstract}
\tableofcontents
\section{Introduction and presentation of the model}

The model of directed polymers in random environment appeared first in the physics literature as
an attempt to modelize roughening in domain wall in the $2D$-Ising model due to impurities \cite{HH}.
It then reached the mathematical community in \cite{IS}, and in \cite{B}, where the author applied
martingale techniques that have became the major technical tools in the study of this model since then.
A lot of progress has been made recently in the mathematical understanding of directed polymer model (see for example \cite{J,CY, CSY, CY_cmp, CH_al, CH_ptrf, CV} and \cite{CSY_rev} for a recent review). 
It is known that there is a phase transition from a delocalized phase at high temperature, where the behavior
of the polymer is diffusive, to a localized phase, where it is expected that the influence of the media is relevant in order to produce nontrivial phenomenons, such as super-diffusivity. These two different situations are
usually referred to as weak and strong disorder, respectively. A simple characterization of this dichotomy is
given in terms of the limit of a certain positive martingale related to the partition function
of this model.

It is known that in low dimensions ($d=1$ or $2$), the polymer is essentially
in the strong disorder phase (see \cite{Lac}, for more precise results), but for $d\geq 3$, there
is a nontrivial region of temperatures where weak disorder holds. A weak form of invariance principle is
proved in \cite{CY}.

However, the exact value of the critical temperature which separates the two regions (when it is finite) remains an open question.
It is known exactly in the case of directed polymers on the tree, where a complete analysis is available
(see \cite{Buf, Fr, KP}). In the case of $\Z^d$, for $d\ge 3$, an $L^2$ computation yields an upper bound on the critical temperature, which is however known not to coincide with this bound (see \cite{BS, Birk} and \cite{CC}). 

We choose to study the same model of directed polymers on diamond hierarchical lattices. These lattices present
a very simple structure allowing to perform a lot of computations together with a richer geometry
than the tree (see Remark \ref{convexity-tree} for more details).   They have been introduced in physics in order to perform exact renormalization group computations for spin systems (\cite{Migd, Kad}). A detailed treatment of more general hierarchical lattices can be found
in \cite{KG1} and \cite{KG2}. For an overview of the extensive literature on Ising and Potts models on hierarchical lattices, we refer the reader to \cite{BZ, DF} and references therein. Whereas statistical mechanics model on trees have to be considered as mean-field versions of the original models, the hierarchical lattice models are in many sense very close to the models on $\Z^d$; they are a very powerful tool to get an intuition for results and proofs on the more complex $\Z^d$ models (for instance, the work on hierarchical pinning model in \cite{GLT} lead to a solution of the original model in \cite{DGLT}. In the same manner, the present work has been a great source of inspiration for \cite{Lac}).

Directed polymers on hierarchical lattices (with bond disorder) appeared in \cite{CD, DG,DGr, DHV}
(see also \cite{R_al} for directed first-passage percolation).
More recently, these lattice models raised the interest of mathematicians in the study of random resistor  networks (\cite{W}), pinning/wetting transitions (\cite{GLT,L}) and diffusion on a percolation cluster (\cite{HK}).

We can also mention \cite{Ham} where the authors consider a random analogue of the hierarchical lattice, where
at each step, each bond transforms either into a series of two bonds or into two bonds in parallel, with probability
$p$ and $p-1$ respectively. 

\vspace{2ex}

Our aim in this paper is to describe the properties of the quenched free energy of directed polymers on
hierarchical lattices with site disorder at high temperature:
\begin{itemize}
	\item First, to be able to decide, in all cases, if the quenched and annealed free energy differ at low temperature.
	\item If they do, we want to be able to describe the phase transition and to compute the critical exponent.
\end{itemize}

We choose to focus on the model with site disorder, whereas \cite{GM_h, CD} focus on the model with \textsl{bond disorder} where computations are simpler.
We do so because we believe that this model is closer to the model of directed polymer in $\Z^{d}$ (in particular, because of the inhomogeneity of the Green Function), and because there exists a nice recursive construction of the partition functions in our case, that leads to a martingale property. Apart from that, both models are very similar, and we will shortly talk about the bound disorder model in section \ref{bddis}.

The diamond hierarchical lattice $D_n$ can be constructed recursively:

\begin{itemize}
	\item $D_0$ is one single edge linking two vertices $A$ and $B$.
	\item $D_{n+1}$ is obtained from $D_n$ by replacing each edges by $b$ branches of $s-1$ edges.
\end{itemize}

\begin{figure}[h]
\begin{center}
\leavevmode
\epsfysize =6.5 cm
\psfragscanon
%\psfrag{wall}[c]{{\tiny The path $\sigma$ (the {\sl wall})}}
\psfrag{a}[c]{\normalsize A}
\psfrag{b}[c]{B}
\psfrag{l0}[c]{$D_0$}
\psfrag{l1}[c]{$D_1$}
\psfrag{l2}[c]{$D_2$}
\psfrag{traject}[l]{{\tiny (a directed path on $D_2$)}}
\epsfbox{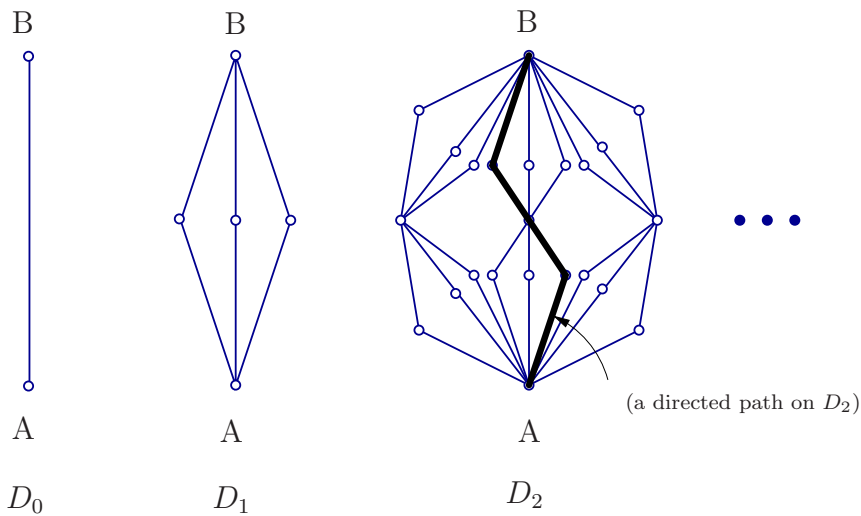}
\end{center}
\caption{\label{fig:Dn} We present here the recursive construction of the first three levels of the hierarchical lattice  $D_n$, for $b=3$, $s=2$.}
\end{figure}
\noindent We can, improperly, consider $D_n$ as a set of vertices, and, with the above construction, we have $D_n\subset D_{n+1}$. We set $D=\bigcup_{n\ge 0} D_n$.
The vertices introduced at the $n$-th iteration are said to belong to the $n$-th generation $V_n=D_n\setminus D_{n-1}$.
We easily see that $|V_n|= (bs)^{n-1} b (s-1)$.

\noindent We restrict to $b\geq 2$ and $s\geq 2$. The case $b=1$ (resp. $s=1$) is not interesting as it just corresponds to a familly of edges
in serie (resp. in parallel)

 We introduce disorder in the system as a set of real numbers associated to vertices $\go=(\go_z)_{z\in D\setminus{\{A,B\}}}$.
Consider $\Gamma_n$ the space of directed paths in $D_n$ linking $A$ to $B$. For each
$g \in \Gamma_n$ (to be understood as a sequence of connected vertices in $D_n$, $(g_0=A, g_1, \dots, g_{s^n}=B)$), 
we define the Hamiltonian

\begin{eqnarray}
	H^{\omega}_n(g) := \sum^{s^{n}-1}_{t=1} \omega(g_t).
	\label{eq: Hn}
\end{eqnarray}

\noindent For $\beta > 0$, $n\geq 1$, we define the (quenched) polymer measure on $\gG_n$ which chooses a path $\gga$ at random  with law

\begin{eqnarray}
	\mu^{\omega}_{\beta, n}(\gamma=g) := \frac{1}{Z_{n}(\beta)} \exp( \beta H^{\omega}_n(g)),
	\label{eq:polymer}
\end{eqnarray}

\noindent where

\begin{eqnarray}
	Z_{n}(\beta)= Z_{n}(\beta, \omega):=\sum_{g \in \Gamma_n} \exp (\beta H^{\omega}_n(g)),
	\label{eq:Zn}
\end{eqnarray}

\noindent is the partition function, and $\gb$ is the inverse temperature parameter.

 In the sequel, we will focus on the case where $\go=( \omega_z,\, z\in D\setminus\{A,B\} )$ is a collection of i.i.d.\ random variables and denote the product measure by $Q$. Let $\go_0$ denote a one dimensional marginal of $Q$, we assume that $\go_0$ has expectation zero, unit variance, and that

\begin{eqnarray}
	\lambda(\beta):= \log Q e^{\beta \go_0}<\infty \quad \forall \gb>0.
	\label{eq:lambda}
\end{eqnarray}
As usual, we define the quenched free energy (see Theorem \ref{thermo}) by

\begin{eqnarray}
	p(\beta) := \lim_{n\to +\infty} \frac{1}{s^{n}} Q \log Z_n(\beta),
	\label{eq:quenched}
\end{eqnarray}

\noindent and its annealed counterpart by

\begin{eqnarray}
	f(\beta) := \lim_{n\to +\infty} \frac{1}{s^{n}} \log Q Z_n(\beta).
	\label{eq:annealed}
\end{eqnarray}

\noindent This annealed free energy can be exactly computed. We will prove

\begin{eqnarray}
	f(\beta):= \lambda(\beta)+ \frac{\log b}{s-1}.
	\label{eq:annealedexact}
\end{eqnarray}

\vspace{3ex}

This model can also be stated as a random dynamical system: given two integer parameters $b$ and $s$ larger than $2$, $\gb>0$, consider the following recursion:

\begin{align}
W_0&\stackrel{\mathcal L}{=}1\notag\\
W_{n+1}&\stackrel{\mathcal L}{=}\frac{1}{b}\sum_{i=1}^b \prod_{j=1}^{s}W_n^{(i,j)}\prod_{i=1}^{s-1}A^{(i,j)}_n\label{eq:rec},
\end{align}

\noindent where equalities hold in distribution, $W_n^{(i,j)}$ are independent copies of $W_n$, and $A^{(i,j)}_n$ are i.i.d.\ random variables, independent of the $W_n^{(i,j)}$ with law

\begin{align*}
A\stackrel{\mathcal L}{=}\exp(\gb\go-\gl(\gb)).
\end{align*}

\noindent In the directed polymer setting, $W_n$ can be interpretative as the normalized partition function

\begin{eqnarray}
	W_n(\beta)=W_n(\beta,\omega) = \frac{Z_n(\beta,\omega)}{Q Z_n(\beta, \omega)}.
	\label{eq:W}
\end{eqnarray}

\noindent Then, (\ref{eq:rec}) turns out to be an almost sure equality if we interpret $W_n^{(i,j)}$ as the partition function of the 
$j$-th edge of the $i$-th branch of $D_1$.

\noindent The sequence $(W_n)_{n\ge 0}$ is a martingale with respect to $\mathcal{F}_n = \sigma(\omega_z:\, z\in \cup^n_{i=1}V_i)$ and
as $W_n>0$ for all $n$, we can define the almost sure limit
$W_{\infty}= \lim_{n\to +\infty} W_n$. Taking limits in both sides of (\ref{eq:rec}), we obtain a functional equation for $W_{\infty}$.

%--------------------------------------------------------------------------------------------------------------------------

%--------------------------------------------------------------------------------------------------------------------------

%--------------------------------------------------------------------------------------------------------------------------

\section{Results}

Our first result is about the existence of the free energy.

\begin{theorem}\label{thermo}
For all $\beta$, the limit

\begin{eqnarray}
	\lim_{n\to +\infty} \frac{1}{s^n} \log Z_n(\beta),
	\label{eq:free}
\end{eqnarray}

\noindent exists a.s. and is a.s. equal to the quenched free
energy $p(\beta)$. In fact for any $\gep>0$, one can find $n_0(\gep,\gb)$ such that
\begin{eqnarray}\label{eq:concentrate}
  Q\left(\left|Z_n-Q\log Z_n\right|> s^n \gep\right)\le \exp\left(-\frac{\gep^{2/3}s^{n/3}}{4}\right), \quad \text{for all } n\ge n_0
\end{eqnarray}

Moreover, $p(\cdot)$ is a strictly convex function of $\beta$.
\end{theorem}

\vspace{4ex}

\begin{rem}\rm
The inequality (\ref{eq:concentrate}) is the exact equivalent of \cite[Proposition 2.5]{CSY}, 
and the proof given there can easily be adapted to our case. It applies concentration results for
martingales from \cite{LV}.
It can be improved in order to obtain the same bound as for Gaussian
environments stated in \cite{CH_ptrf} (see \cite{Cnotes} for details). However, it is believed that it is
not of the optimal order, similar to the case of directed polymers on $\Z^d$.
\end{rem}

\begin{rem}\rm \label{convexity-tree}
The strict convexity of the free energy is an interesting property. It is known that it holds also for the directed polymer on $\Z^d$ but not on the tree. In the 
later case, the free energy is strictly convex only for values of $\beta$ smaller than the critical value $\beta_c$ 
(to be defined latter) and it is linear on $[\beta_c,+\infty)$. This fact is related to the particular structure of
the tree that leads to major simplifications in the 'correlation' structure of the model (see \cite{Buf}). The strict convexity, in
our setting, arises essentially from the property that two path on the hierarchical lattice can re-interesect after being separated at some step. This underlines once more, that $\Z^d$ and the hierarchical lattice have a lot of features in common, which they do not share with the tree.
\end{rem}

We next establish the martingale property for $W_n$ and the zero-one law for its limit.

\begin{lemma}\label{martingale}
$(W_n)_n$ is a positive $\mathcal{F}_n$-martingale. It converges
$Q$-almost surely to a non-negative limit $W_{\infty}$ that satisfies
the following zero-one law:

\begin{eqnarray}
 Q \left( W_{\infty} > 0 \right)\,\in\, \{0,1\}.
 \label{zeroone}
\end{eqnarray}
\end{lemma}

\vspace{2ex}
\noindent Recall that martingales appear when the disorder is displayed on sites, in contrast with disorder on bonds 
as in \cite{CD,DG}.
\vspace{3ex}

Observe that

\begin{eqnarray*}
p(\beta)-f(\beta)=\lim_{n\to +\infty} \frac{1}{s^n} \log W_n(\beta),
\end{eqnarray*}

\noindent so, if we are in the situation $Q(W_{\infty}>0)=1$, we have
that $p(\beta)=f(\beta)$. This motivates the following definition:

\begin{definition}\label{disorder}
 If $Q(W_{\infty}>0)=1$, we say that weak disorder holds. In the opposite situation,
 we say that strong disorder holds.
\end{definition}

\begin{rem}\rm
Later, we will give a statement(Proposition \ref{th:fracmom}) that guarantees that strong disorder is equivalent to $p(\gb)\neq f(\gb)$, a situation that is sometimes called very strong disorder. This is believed to be true for polymer models on $\Z^d$ or $\R^d$ but it remains an unproved and challenging conjecture in dimension $d\ge 3$ (see \cite{CH_al}).
\end{rem}

The next proposition lists a series of partial results that in some sense clarify
the phase diagram of our model. 

\begin{proposition}\label{misc}

  $(i)$ There exists $\beta_0\in[0, +\infty]$ such that strong disorder holds
  for $\beta > \beta_0$ and weak disorder holds for $\beta\le  \beta_0$.
  \vspace{2ex}

  $(ii)$ If $b>s$, $\beta_0>0$. Indeed, there exists $\beta_2\in(0,\infty]$ such that for all $\beta<\beta_2$,
         $\sup_n Q(W^2_{n}(\beta))<+\infty$, and therefore weak disorder holds.
  \vspace{2ex}

  $(iii)$ If $\beta \lambda'(\beta)-\lambda(\beta)> \frac{2\log b}{s-1}$, then strong
          disorder holds.
  \vspace{2ex}

  $(iv)$ %If $b+1<2s$, then $\beta_0 > \beta_2$.
  In the case where $\go_z$ are gaussian random variables, $(iii)$ can be improved for $b>s$: strong disorder holds as soon as $\gb>\sqrt{\frac{2(b-s)\log b}{(b-1)(s-1)}}$.
  \vspace{2ex}

  $(v)$ If $b\leq s$, then strong disorder holds for all $\beta$.
\end{proposition}

\begin{rem}\rm
On can check that the formula in $(iii)$ ensures that $\gb_0<\infty$ whenever the distribution of $\go_z$ is unbounded.
\end{rem}
\begin{rem}\rm
An implicit formula is given for $\gb_2$ in the proof and this gives a lower bound for $\gb_0$. However, when $\gb_2<\infty$, it never coincides with the upper bound given by $(iii)$ and $(iv)$, and therefore knowing the exact value of the critical temperature when $b>s$ remains an open problem.
\end{rem}
\vspace{4ex}

We now provide more quantitative information for the regime considered in $(v)$:

\begin{theorem} \label{th:bs}
When $s>b$, there exists a constant $c_{s,b}=c$ such that for any $\gb\le 1$ we have
\begin{align*}
\frac{1}{c}\gb^{\frac{2}{\alpha}}\le \lambda(\beta)-p(\beta) \le c\gb^{\frac{2}{\alpha}}
\end{align*}
where $\alpha=\frac{\log s-\log b}{\log s}$.
\end{theorem}

\begin{theorem} \label{th:ss}
When $s=b$, there exists a constant $c_s=c$ such that for any $\gb\le 1$ we have
\begin{align*}
\exp\left(-\frac{c}{\gb^2}\right)\le  \lambda(\beta) - p(\beta) \le c\exp\left(-\frac{1}{c\gb}\right)
\end{align*}
\end{theorem}

\vspace{4ex}
In the theory of directed polymer in random environment, it is believed that, in low dimension, the quantity $\log Z_n$ undergoes large fluctuations around its average (as opposed to what happens in the weak disorder regime where the fluctuation are of order $1$). More precisely: it is believed that there exists  exponents $\xi>0$ and $\chi\ge 0$ such that
\begin{equation}
 \log Z_n-Q\log Z_n \asymp N^{\xi} \text{ and } \var_Q  \log Z_n\asymp N^{2\chi},
\end{equation}
where $N$ is the length of the system ($=n$ on $\Z^d$ and $s^n$ one our hierarchical lattice).
In the non-hierarchical model this exponent is of major importance as it is closely related to the {\sl volume exponent} $\xi$ that gives
the spatial fluctuation of the polymer chain (see e.g.\ \cite{J} for a discussion on fluctuation exponents). Indeed it is conjectured for the $ \Z^d$ models that
\begin{equation}
\chi=2\xi-1.
\end{equation}
This implies that the polymer trajectories are superdiffusive as soon as $\chi>0$.
  In our hierarchical setup, there is no such geometric interpretation but having a lower bound on the fluctuation allows to get a significant localization result.

\begin{proposition}\label{fluctuations}
When $b<s$, there exists a constant $c$ such that for all $n\ge 0$ we have
\begin{equation}
\var_Q\left(\log Z_n\right)\ge \frac{c (s/b)^{n}}{\gb^2}.
\end{equation}
Moreover, for any $\gep>0$, $n\ge 0$, and $a\in \bbR$,
\begin{equation}
Q\left\{\log Z_n \in [a, a+\gep(s/b)^{n/2}]\right\}\le \frac{8\gep}{\gb}.
\end{equation}
\end{proposition}

\noindent This implies that if the fluctuation exponent $\chi$ exists, $\chi\ge \frac{\log s -\log b}{2\log s}$.
We also have the corresponding result for the case $b=s$

\begin{proposition}\label{fluc2}
When $b=s$, there exists a constant $c$ such that for all $n\ge 0$ we have
\begin{equation}
\var_Q \left(\log Z_n\right)\ge \frac{c n}{\gb^2}.
\end{equation}
Moreover for any $\gep>0$, $n\ge 0$, and $a\in \bbR$,
\begin{equation}
Q\left\{\log Z_n \in [a, a+\gep\sqrt{n}]\right\}\le \frac{8\gep}{\gb}.
\end{equation}
\end{proposition}

From the fluctuations of the free energy we can prove the following:
For $g \in \gG_n$ and $m<n$, we define $g|_m$ to be the restriction of $g$ to $D_m$.

\begin{cor}\label{locloc}
If $b\le s$, and $n$ is fixed we have
\begin{equation}
 \lim_{n\to\infty} \sup_{g\in \gG_m} \mu_n(\gga|_m=g)=1,
\end{equation}
\noindent where the convergence holds in probability.
\end{cor}
Intuitively this result means that if one look on a large scale, the law of $\mu_n$ is concentrated in the neighborhood of a single path. Equipping $\gG_n$ with a natural metric (two path $g$ and $g'$ in $\gG_n$ are at distance $2^{-m}$ if and only if $g|_m\ne g'|_m$ and $g|_{m-1}=g|_{m-1}$) makes this statement rigorous.

\begin{rem}\rm
Proposition \ref{misc}$(v)$ brings the idea that $b\le s$ for this hierarchical model is equivalent to the $d\le 2$ case for the model in $\Z^{d}$ (and that $b>s$ is equivalent to $d>2$). Let us push further the analogy: let $\gga^{(1)}$ , $\gga^{(2)}$ be two paths chosen uniformly at random in $\gG_n$ (denote the uniform-product law by $P^{\otimes 2}$), their expected site overlap is of order $(s/b)^n$ if $b<s$, of order $n$ if $b=s$, and of order $1$ if $b>s$.
If one denotes by $N= s^n$ the length of the system, one has
\begin{equation}
 P^{\otimes 2}\left[\sum_{t=0}^{N} \ind_{\{\gamma^{(1)}_t= \gamma^{(2)}_t\}}\right]\, \asymp\, \begin{cases}N^{\alpha} \quad \text{ if } b < s,\\
                                                          \log N \quad \text{ if } b=s,\\
								1 \text{ if } b> s,
                                                          \end{cases}
 \end{equation}
(where $\alpha=(\log s-\log b)/\log s$).
Comparing this to the case of random walk on $\Z^d$, we can infer that the case $b=s$ is just like $d=2$ and that the case $d=1$ is similar to $b=\sqrt{s}$ ($\ga=1/2$). One can check in comparing \cite[Theorem 1.4, 1.5, 1.6]{Lac} with Theorem \ref{th:bs} and \ref{th:ss}, that this analogy is relevant.
\end{rem}

The paper is organised as follow
\begin{itemize}
 \item In section \ref{mtricks} we prove %Theorem \ref{thermo}
some basic statements about the free energy, Lemma \ref{martingale} and
the first part of Proposition \ref{misc}.
 \item Item $(ii)$ from Proposition \ref{misc} is proved in Section $5.1$.
Item $(v)$ is a consequence of Theorems \ref{th:bs} and \ref{th:ss}.
 \item Items $(iii)$ and $(iv)$ are proved in Section $6.3$.
Theorems \ref{th:bs} and \ref{th:ss} are proved in Section $6.1$ and $6.3$
respectively.
 \item In section \ref{flucloc} we prove Propositions \ref{fluctuations} and \ref{fluc2} and Corrolary \ref{locloc}.
 \item  In section \ref{weakpol} we define and investigate the properties of the infinite volume polymer measure in the weak disorder phase.
 \item In section \ref{bddis} we shortly discuss about the bond disorder model.
\end{itemize}

%--------------------------------------------------------------------------------------------------------------------------

%--------------------------------------------------------------------------------------------------------------------------

%--------------------------------------------------------------------------------------------------------------------------

\section{Martingale tricks and free energy}\label{mtricks}

We first look at to the existence of the quenched free energy

\begin{eqnarray*}
  p(\beta)= \lim_{n\to +\infty} \frac{1}{n} Q \left( \log Z_n(\beta) \right),
\end{eqnarray*}

\noindent and its relation with the annealed free energy. The case $\beta=0$ is
somehow instructive. It gives the number of paths in $\Gamma_n$ and is
handled by the simple recursion:

\begin{eqnarray*}
Z_n(0)=b \left( Z_{n-1}(0) \right).
\end{eqnarray*}

\noindent This easily yields

\begin{eqnarray}
\label{paths}
 |\Gamma_n| = Z_n(\beta=0)= b^{\frac{s^{n}-1}{s-1}}.
\end{eqnarray}

Much in the same spirit than (\ref{eq:rec}), we can find a recursion for
$Z_n$:

\begin{eqnarray}
	Z_{n+1}= \sum^b_{i=1} Z^{(i,1)}_n \cdots Z^{(i,s)}_n \times e^{\beta \omega_{i,1}}\cdots e^{\beta \omega_{i,s-1}}.
	\label{eq:recz}
\end{eqnarray}
The existence of the quenched free energy follows by monotonicity: we have

$$Z_{n+1} \geq Z^{(1,1)}_n Z^{(1,2)}_n \cdots Z^{(1,s)}_n \times e^{\beta \omega_{1,1}}\cdots e^{\beta \omega_{1,s-1}},$$

\noindent so that (recall the $\omega$'s are centered random variables)

$$\frac{1}{s^{n+1}} Q \log Z_{n+1} \geq \frac{1}{s^n} Q \log Z_n.$$

\noindent The annealed free energy provides an upper bound:

\begin{eqnarray*}
\frac{1}{s^n} Q \log Z_n &\leq& \frac{1}{s^n} \log Q Z_n\\
                         &=& \frac{1}{s^n} \log e^{\lambda(\beta)(s^n-1)}Z_n(\beta=0)\\
                         &=& \left( 1 - \frac{1}{s^n} \right) \left( \lambda(\beta) + \frac{\log b}{s-1} \right)\\
                         &=& \left( 1 - \frac{1}{s^n} \right) f(\beta).
\end{eqnarray*}

\vspace{3ex}

%------------------------------------------------CONVEXITY OF THE FREE ENERGY-

We now prove the strict convexity of the free energy. The proof
is essentially borrowed from \cite{CPV}, but it is remarkably simpler
in our case.

\begin{proof}[Proof of the strict convexity of the free energy]
We will consider a Bernoulli environment ($\go_z=\pm 1$ with probability $p$, $1-p$; note that our assumptions on the variance and expectation for $\go$ are violated but centering and rescaling $\go$ does not change the argument).
We refer to \cite{CPV} for generalization to more general environment.

\noindent An easy computation yields

\begin{eqnarray*}
\frac{d^2}{d\beta^2} Q \log Z_n = Q {\rm Var}_{\mu_n} H_n (\gamma).
\end{eqnarray*}

\noindent We will prove that for each $K>0$, there exists a constant
$C$ such that, for all $\beta \in[0,K]$ and $n\geq 1$,

\begin{eqnarray}\label{lowerboundenergy}
{\rm Var}_{\mu_n} H_n (\gamma) \geq C s^n
\end{eqnarray}

\vspace{2ex}
For $g \in \gG_n$ and $m<n$, we define $g|_m$ to be the restriction of $g$ to $D_m$.
By the conditional variance formula,

\begin{eqnarray}
\nonumber
{\rm Var}_{\mu_n} H_n &=& \mu_n \left( {\rm Var}_{\mu_n} (H_n(\gamma)\, |\, \gamma_{|_{n-1}}) \right) +
                   {\rm Var}_{\mu_n} \left( \mu_n(H_n(\gamma)\, |\,  \gamma_{|_{n-1}})\right)\\
                   \label{condvar}
                  &\geq& \mu_n \left( {\rm Var}_{\mu_n} (H_n(\gamma) \,|\, \gamma_{|_{n-1}}) \right)
\end{eqnarray}

\noindent Now, for $l=0,...,s^{n-1}-1$, $g \in \Gamma_n$, define

\begin{eqnarray*}
H^{(l)}_n (g) = \sum^{(l+1)s-1}_{t=ls+1} \omega(g_t),
\end{eqnarray*}

\noindent so  (\ref{condvar}) is equal to

\begin{eqnarray*}
\mu_n  {\rm Var}_{\mu_n} \left( \sum^{s^{n-1}-1}_{l=0} H^{(l)}_{n}(\gamma) | \gamma_{|_{n-1}} \right)
= \sum^{s^{n-1}-1}_{l=0} \mu_n  {\rm Var}_{\mu_n} \left(  H^{(l)}_{n}(\gamma) | \gamma_{|_{n-1}} \right),
\end{eqnarray*}

\noindent by independence. Summarizing,

\begin{equation}\label{varcond}
{\rm Var}_{\mu_n} H_n\geq \sum^{s^{n-1}}_{l=1} \mu_n  {\rm Var}_{\mu_n} \left(  H^{(l)}_{n}(\gamma) | \gamma_{|_{n-1}} \right).
\end{equation}

\noindent The rest of the proof consists in showing that each term of the sum is bounded from below by a positive constant,
uniformly in $l$ and $n$.
For any $x\in D_{n-1}$ such that the graph distance between $x$ and $A$ is $ls$ in $D_n$ (i.e.\ $x\in D_{n-1}$), we define the set of environment

\begin{eqnarray*}
M(n,l,x)= \left\lbrace \omega: \left|\{H^{(l)}_n(g,\omega)\,:\, g\in \gG_n, g_{ls}=x\} \right| \geq 2 \right\rbrace.
\end{eqnarray*}

\noindent These environments provide the fluctuations in the energy needed for the
uniform lower bound we are searching for.
One second suffices to convince oneself
that $Q(M(n,l,x)>0$, and does not depend on the parameters $n,\, l$ or $x$. Let $Q(M)$ denote improperly the common value of $Q( M(n,l,x))$. 
Now, it is easy to see (from \eqref{varcond}) that there exists a constant $C$ such that for all $\gb<K$,

\begin{eqnarray*}
Q \left[{\rm Var}_{\mu_n} H_n\right] &\ge& C Q \left[\sum^{s^{n-1}-1}_{l=1} \sum_{x\in D_{n-1}}  {\bf 1}_{M(n, l, x)}\mu_n ({\gamma_{ls} = x})\right].
    \end{eqnarray*}
    
\noindent Define now $\mu^{(l)}_n$ as the polymer measure in the environment obtained from
$\omega$ by setting $\omega(y)=0$ for all sites $y$ which distance to $0$ is between $ls$ and $(l+1)s$.
One can check that for all $n$, and all path $g$,

\begin{eqnarray*}
\exp(-2\gb(s-1)) \mu^{(l)}_n(\gamma=g) \leq \mu_n(\gamma=g) \leq \exp(2\gb (s-1)) \mu^{(l)}_n(\gamma).
\end{eqnarray*}

\noindent We note that under $Q$, $\mu_n^{(l)}(\gamma_{ls} = x)$ and $\ind_{M(n, l, x)}$ are random variables, so that 

\begin{eqnarray*}
Q\left[ {\rm Var}_{\mu_n} H_n \right]&\geq& C\exp(-2\gb(s-1)) Q \left[\sum^{s^{n-1}-1}_{l=0} \sum_x  {\bf 1}_{M(n, l, x)}\mu^{(l)}_n ({\gamma_{ls} = x})\right]\\
   &=& C\exp(-2\gb(s-1))\sum^{s^{n-1}}_{l=1} \sum_{x\in D_{n-1}} Q(M(n,l,x))   Q\left[ \mu^{(l)}_n ({\gamma_l = x})\right]\\
   &=& C\exp(-2\gb(s-1)) Q(M) s^{n-1}.
\end{eqnarray*}
\end{proof}

\vspace{3ex}

%----------------------------------------------------MARTINGALE PROPERTY-

\noindent We now establish the martingale property for the normalized free energy.

\begin{proof}[Proof of Lemma \ref{martingale}]

Set $z_n=Z_n(\beta=0)$. We have already remarked that this is just the number of (directed) paths in $D_n$, 
and its value is given by (\ref{paths}). 
Observe that
$g\in \gG_n$ visits $s^n(s-1)$ sites of $n+1$-th generation. The restriction of paths in $D_{n+1}$
to $D_n$ is obviously not one-to-one as for each path $g'\in \gG_n$, there are $b^{s^n}$ paths
in $\gG_{n+1}$ such that $g|_n=g'$. Now,

\begin{eqnarray*}
Q\left(Z_{n+1}(\beta)| \mathcal{F}_n\right)
&=& \sum_{g \in D_{n+1}} Q\left( e^{\beta H_{n+1}(g)}|\mathcal{F}_n\right)\\
&=&\sum_{g' \in D_n} \sum_{g \in D_{n+1}}Q\left( e^{\beta H_{n+1}(g)}|\mathcal{F}_n\right){\bf 1}_{g|_n=g'}\\
&=&\sum_{g' \in D_n}\sum_{g \in D_{n+1}} e^{\beta H_{n}(g')}e^{s^n(s-1)\lambda(\beta)}{\bf 1}_{g|_n=g'}\\
&=&\sum_{g' \in D_n} e^{\beta H_{n}(g')}e^{s^n(s-1)\lambda(\beta)}\sum_{g \in D_{n+1}}{\bf 1}_{g|_n=g'}\\
&=&e^{s^n(s-1)\lambda(\beta)}b^{s^n}\sum_{g' \in D_n} e^{\beta H_{n}(g')}\\
&=&Z_n(\beta) \frac{z_{n+1}e^{s^{n+1}\lambda(\beta)}}{z_n e^{s^n \lambda(\beta)}}.
\end{eqnarray*}

\noindent This proves the martingale property. For (\ref{zeroone}), let's generalize a little
the preceding restriction procedure. As before, for a path $g \in D_{n+k}$, denote by
$g |_n$ its restriction to $D_n$. Denote by $I_{n,n+k}$ the set of time indexes that have been
removed in order to perform this restriction and by $N_{n,n+k}$ its cardinality. Then

\begin{eqnarray*}
Z_{n+k}= \sum_{g \in D_n} e^{\beta H_n(g)} \sum_{g' \in D_{n+k}, g'|_n=g}
\exp \left\lbrace \beta \sum_{t\in I_{n,n+k}} \omega(g'_t) \right\rbrace.
\end{eqnarray*}

\noindent Consider the following notation, for $g\in \gG_n$,
\begin{eqnarray*}
\tilde{W}_{n,n+k}(g)= c^{-1}_{n,n+k}\sum_{g' \in D_{n+k}, g'|_n=g}\exp \left\lbrace \beta \sum_{t\in I_{n,n+k}} \omega(g'_t)  - N_{n,n+k} \lambda(\beta)\right\rbrace,
\end{eqnarray*}

\noindent where $c_{n,n+k}$ stands for the number paths in the sum. With this notations, we have, 

\begin{eqnarray}
 W_{n+k} = \frac{1}{z_n} \sum_{g \in D_n} e^{\beta H_n(g)-(s^{n}-1)\lambda(\beta)}\tilde{W}_{n,n+k}(g),
\end{eqnarray}

\noindent and, for all $n$,

\begin{eqnarray}
\left\lbrace W_{\infty}=0 \right\rbrace =
\left\lbrace \tilde{W}_{n,n+k}(g) \to 0,\, {\rm as}\, k\to +\infty,\, \forall \, g \in D_n \right\rbrace.
\label{eq:WW}
\end{eqnarray}

\noindent The event in the right hand side is measurable with respect to the disorder of generation
not earlier than $n$. As $n$ is arbitrary, the right hand side of (\ref{eq:WW}) is in the tail $\sigma$-algebra
and its probability is either $0$ or $1$.
\end{proof}

This, combined with FKG-type arguments (see \cite[Theorem 3.2]{CY} for details), proves part $(i)$ of
Proposition \ref{misc}. Roughly speaking, the FKG inequality is used to insure that there is no reentrance
phase.

%--------------------------------------------------------------------------------------------------------------------------

%--------------------------------------------------------------------------------------------------------------------------

%--------------------------------------------------------------------------------------------------------------------------

%--------------------------------------------------------------------------------------------------------------------------

%--------------------------------------------------------------------------------------------------------------------------

%--------------------------------------------------------------------------------------------------------------------------

\section{Second moment method and lower bounds}

This section contains all the proofs concerning coincidence of annealed and quenched free--energy for $s>b$ and lower bounds on the free--energy for $b\le s$ (i.e. half of the results from Proposition \ref{misc} to Theorem \ref{th:ss}.)
First, we discuss briefly the condition on $\gb$ that one has to fulfill to to have $W_{n}$ bounded in $\bbL_2(Q)$. Then for the cases when strong disorder holds at all temperature ($b\le s$), we present a method that combines control of the second moment up to some scale $n$ and a percolation argument to get a lower bound on the free energy.
\\
First we investigate how to get the variance of $W_n$ (under $Q$).
From \eqref{eq:rec} we get the induction for the variance $v_n=Q\left[(W_n-1)^2\right]$:

\begin{eqnarray}
v_{n+1}&=&\frac{1}{b}\left(e^{(s-1)\gga(\gb)}(v_n+1)^s-1\right), \label{eq:var}\\
v_0&=&0.\label{eq:v00}
\end{eqnarray}
where $\gga(\gb):=\gl(2\gb)-2\gl(\gb)$.

\subsection{The $L^2$ domain: $s<b$}

If $b>s$, and $\gga(\gb)$ is small, the map
\begin{align*}
g:\ x\mapsto\frac{1}{b}\left(e^{(s-1)\gga(\gb)}(x+1)^s-1\right)
\end{align*}
possesses a fixed point. In this case, \eqref{eq:var} guaranties that $v_n$ converges to some finite limit. Therefore, in this case, $W_n$ is a positive martingale bounded in $\bbL^2$, and therefore converges almost surely to $W_\infty \in \bbL^2(Q)$ with $Q W_{\infty}=1$, so that
\begin{align*}
p(\beta)-\lambda(\beta)=\lim_{n\rightarrow\infty}\frac{1}{s^n}\log W_n=0,
\end{align*}

\noindent and weak disorder holds.
One can check that $g$ has a fixed point if and only if
\begin{align*}
\gga(\gb)\le \frac{s}{s-1}\log \frac{s}{b}-\log \frac{b-1}{s-1}
\end{align*}

%--------------------------------------------------------------------------------------------------------------------------

%--------------------------------------------------------------------------------------------------------------------------

%--------------------------------------------------------------------------------------------------------------------------

\subsection{Control of the variance: $s>b$}

For $\epsilon>0$, let $n_0$ be the smallest integer such that $v_n\ge \gep$.

\vspace{2ex}

\begin{lemma}\label{th:bss}
For any $\gep>0$, there exists a constant $c_{\gep}$ such that for any $\gb\le 1$
\begin{align*}
n_0\ge \frac{2 |\log \gb|}{\log s - \log b}-c_{\gep}.
\end{align*}
\end{lemma}

\begin{proof}

Expanding \eqref{eq:var} around $\gb=0$, $v_n=0$, we find a constant $c_1$ such that, whenever $v_n\le 1$ and $\gb\le 1$,
\begin{align}
v_{n+1}\le \frac{s}{b}(v_n+c_1\gb^2)(1+c_1v_n) \label{eq:varmod}.
\end{align}
 Using \eqref{eq:varmod}, we obtain by induction
\begin{align*}
v_{n_0}\le \prod_{i=0}^{n_0-1}(1+c_1v_i)\left[c_1\gb^2\left(\sum_{i=0}^{n_0-1}(s/b)^{i}\right)\right].
\end{align*}
From \eqref{eq:var}, we see that $v_{i+1}\ge (s/b)v_i$. By definition of $n_0$, $v_{n_0-1}<\epsilon,$ so that $v_{i} <\gep (s/b)^{i-n_0+1}$. Then
\begin{align*}
\prod_{i=0}^{n_0-1}(1+c_1v_i)\le \prod_{i=0}^{n_0-1}(1+c_1\gep(s/b)^{i-n_0+1})\le \prod_{k=0}^{\infty}(1+c_1\gep(s/b)^{-k})\le 2,
\end{align*}
where the last inequality holds for $\gep$ small enough.
In that case we have
\begin{align*}
\gep\le v_{n_0}\le 2 c_1 \gb^2 (s/b)^{n_0},
\end{align*}
so that \begin{equation*}
         n_0\ge \frac{\log(\gep/2c_1\gb^2)}{\log(s/b)}.
        \end{equation*}
\end{proof}

%--------------------------------------------------------------------------------------------------------------------------

%--------------------------------------------------------------------------------------------------------------------------

%--------------------------------------------------------------------------------------------------------------------------

\subsection{Control of the variance: $s=b$}

\begin{lemma}\label{th:sss}
There exists a constant $c_2$ such that, for every $\gb\le 1$, 
\begin{align*}
v_n\le \gb, \quad \forall \, n\ \le  \ \frac{c_2}{\gb}.
\end{align*}
\end{lemma}
\begin{proof}
%From (\ref{eq:varmod}), we have
%\begin{align*}
%v_{n+1}\le (v_n + c_1\gb^2)(1+c_1 v_n),
%\end{align*}
%for $\beta \leq 1$ and $v_n\leq 1$. 

By (\ref{eq:varmod}) and induction we have, for any $n$ such that $v_{n-1}\le 1$ and $\beta\leq 1$,
\begin{align*}
v_{n}\le n\gb^2 \prod_{i=0}^{n-1} (1+c_1v_i).
\end{align*}
Let $n_0$ be the smallest integer such that $v_{n_0}>\gb$. By the above formula, we have
\begin{align*}
v_{n_0}\le n_0\gb^2(1+c_1\gb)^{n_0}
\end{align*}
Suppose that $n_0\le (c_2/\gb)$, then
\begin{align*}
\gb\le v_{n_0}\le c_2c_1\gb(1+c_1\gb)^{c_2/\gb}.
\end{align*}
If $c_4$ is chosen small enough, this is impossible.
\end{proof}

%--------------------------------------------------------------------------------------------------------------------------

%--------------------------------------------------------------------------------------------------------------------------

%--------------------------------------------------------------------------------------------------------------------------

\subsection{Directed percolation on $D_n$} For technical reasons, we need to get some understanding on directed independent bond percolation on $D_n$.
Let $p$ be the probability that an edge is open (more detailed considerations about edge disorder
are given in the last section). 
The probability of having an open path from $A$ to $B$ in $D_n$ follows the recursion
\begin{align*}
p_0&=p,\\
p_n&=1-(1-p_{n-1}^s)^b.
\end{align*}
On can check that the map $x\mapsto 1-(1-x^s)^b$ has a unique unstable fixed point on $(0,1)$; we call it $p_c$.
Therefore if $p>p_c$, with a probability tending to $1$, there will be an open path linking $A$ and $B$ in $D_n$. If $p<p_c$, $A$ and $B$ will be disconnected in $D_n$ with probability tending to $1$. If $p=p_c$, the probability that $A$ and $B$ are linked in $D_n$ by an open path is stationary. See \cite{HK} for a deep investigation of percolation on hierarchical lattices.

\subsection{From control of the variance to lower bounds on the free energy}

Given $b$ and $s$, let $p_c=p_c(b,s)$ be the critical parameter for directed bond percolation.

\begin{proposition}\label{th:perco}
Let $n$ be an integer such that $v_n=Q (W_n-1)^2 <\frac{1-p_c}{4}$ and $\gb$ such that $p(\gb)\le (1-\log2)$.
Then
\begin{align*}
\lambda(\beta)-p(\beta)\geq s^{-n}
\end{align*}
\end{proposition}
\begin{proof}

If $n$ is such that $Q\left[ (W_n-1)^2\right] <\frac{1-p_c}{4}$, we apply Chebycheff inequality to see that
\begin{align*}
Q(W_n<1/2)\le 4 v_n< 1-p_c.
\end{align*}

Now let be $m\ge n$. $D_m$ can be seen as the graph $D_{m-n}$ where the edges have been replaced by i.i.d.\ copies of $D_n$ with its environment (see fig. \ref{perco}).
To each copy of $D_n$ we associate its renormalized partition function; therefore, to each edge $e$ of $D_{m-n}$ corresponds an independent copy of $W_n$, $W_n^{(e)}$. By percolation (see fig. \ref{perco2}), we will have, with a positive probability not depending on $n$, a path in $D_{m-n}$ linking $A$ to $B$, going only through edges which associated $W_n^{(e)}$ is larger than $1/2$.

\begin{figure}[h]
\begin{center}
\leavevmode
\epsfysize =6.5 cm
\psfragscanon
\psfrag{dn}[c]{\tiny{Dn}}
\psfrag{a}[c]{A}
\psfrag{b}[c]{B}
\psfrag{copies}[c]{\tiny{Independent copies of system of rank $n$.}}
\epsfbox{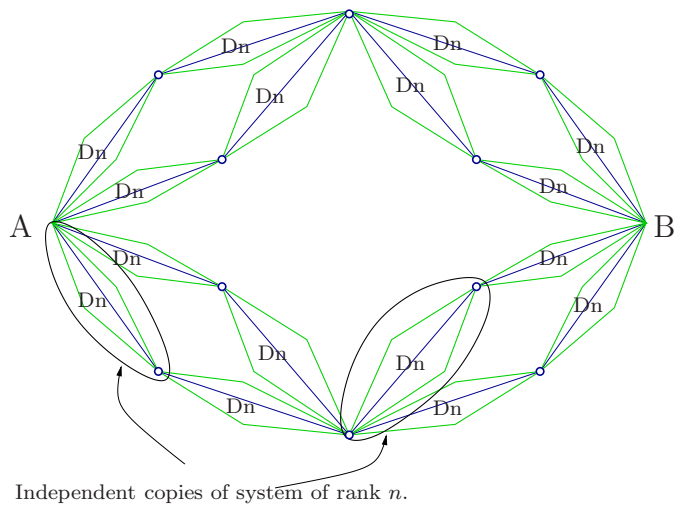}
\end{center}
\caption{\label{perco}On this figure, we scheme how $D_{n+m}$ with its random environment can be seen as independent copies of $D_n$ arrayed as $D_m$.
Here, we have $b=s=2$ $m=2$, each diamond corresponds to a copy of $D_n$ (we can identify it with an edge and get the underlying graph $D_2$). Note that we also have to take into account the environment present on the vertices denoted by circles.}
\end{figure}

\begin{figure}[h]
\begin{center}
\leavevmode
\epsfysize =6.5 cm
\psfragscanon
\psfrag{a}[c]{A}
\psfrag{b}[c]{B}
\psfrag{perco}[c]{\tiny{an open percolation path}}
\epsfbox{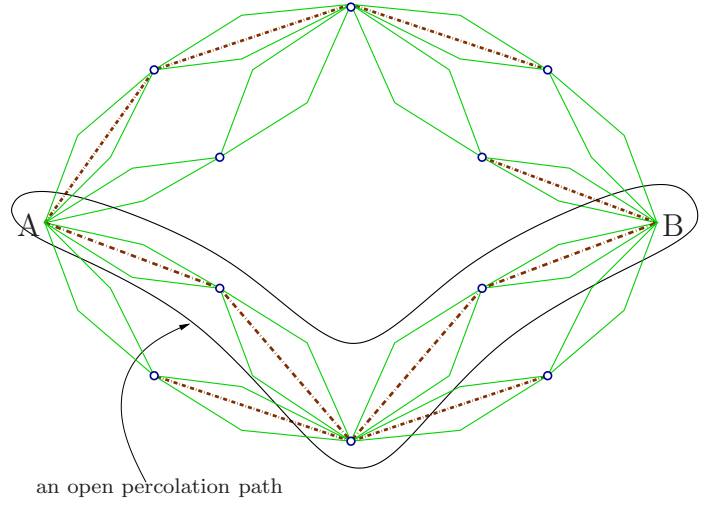}
\end{center}
\caption{\label{perco2}We represent here the percolation argument we use. In the previous figure, we have replaced by an open edge any of the copies of $D_n$ for which satisfies $W_n\ge 1/2$. As it happens with probability larger than $p_c$, it is likely that we can find an open path linking A to B in $D_{n+m}$, especially if $m$ is large.}
\end{figure}

When such paths exist, let $\gamma_0$ be one of them (chosen in a deterministic manner, e.g.\ the lowest such path for some geometric representation of $D_n$).
We look at the contribution of these family of paths in $D_m$ to the partition function.
We have
\begin{align*}
W_m\ge (1/2)^{s^{m-n}}\exp\left(\sum_{z\in \gamma_0} \gb\go_z-\gl(\gb)\right)
\end{align*}
Again, with positive probability (say larger than $1/3$), we have $\sum_{z\in \gamma_0} \go_z\ge 0$ (this can be achieved the the central limit theorem). Therefore with positive probability we have
\begin{align*}
\frac{1}{s^m}\log W_m\ge -\frac{1}{s^{n}}(\log 2+\gl(\gb)).
\end{align*}
As $1/s^m \log W_m$ converges in probability to the free energy this proves the result.

\end{proof}

\begin{proof}[Proof of the right-inequality in Theorems \ref{th:bs} and \ref{th:ss}].

The results now follow by combining Lemma \ref{th:bss} or \ref{th:sss} for $\gb$ small enough, with Proposition
\ref{th:perco}.

\end{proof}

%--------------------------------------------------------------------------------------------------------------------------

%--------------------------------------------------------------------------------------------------------------------------

%----------------------------------------------------------------------------------------------FRACTIONAL MOMENT METHOD---

\section{Fractional moment method, upper bounds and strong disorder}\label{fm}

In this section we develop a way to find an upper bound for $\lambda(\beta)-p(\beta)$, or just to find out if strong disorder hold. The main tool we use are fractional moment estimates and measure changes.

\subsection{Fractional moment estimate}
In the sequel we will use the following notation. Given a fixed parameter $\theta\in(0,1)$, define
\begin{eqnarray}
 u_n&:=&Q W_n^{\theta},\\ \label{un}
 a_\theta&:=&Q A^{\theta}=\exp(\gl(\theta\gb)-\theta\gl(\gb)) \label{aomega}.
\end{eqnarray}

\begin{proposition}\label{th:fracmom}
The sequence $(f_n)_n$ defined by 
$$f_n:=\theta^{-1}s^{-n}\log \left(a_\theta b^{\frac{1-\theta}{s-1}}u_n\right)$$
\noindent is decreasing and we have
\begin{align*}
\lim_{n\rightarrow\infty} f_n\ge p(\beta)-\lambda(\beta).
\end{align*}
\vspace{1ex}

\noindent $(i)$ In particular, if for some $n\in\N$, $u_n< a_{\theta}^{-1}b^{\frac{\theta-1}{s-1}}$, strong disorder holds.
\vspace{1ex}

\noindent $(ii)$ Strong disorder holds in particular if $a_\theta< b^{\frac{\theta-1}{s-1}}$.

%\begin{itemize}
%	\item In particular, if for some $n\in\N$, $u_n< a^{-1}b^{\frac{\theta-1}{s-1}}$, strong disorder holds.
 % \item Strong disorder holds also if $a_\theta< b^{\frac{\theta-1}{s-1}}$.
%\end{itemize}
\end{proposition}

\begin{proof}
The inequality $\left(\sum a_i\right)^{\theta}\le \sum a_i^{\theta}$ (which holds for any $\theta\in (0,1)$ and any collection of positive numbers $a_i$) applied to
\eqref{eq:rec} and averaging with respect to $Q$ gives

\begin{align*}
u_{n+1}\le b^{1-\theta}u_n^{s}a_{\theta}^{s-1}
\end{align*}

\noindent From this we deduce that the sequence
\begin{align*}
s^{-n}\log\left( a_{\theta} b^{\frac{1-\theta}{s-1}}u_n\right)
\end{align*}

\noindent is decreasing. Moreover we have
\begin{align*}
p(\beta)-\lambda(\beta)=\lim_{n\rightarrow\infty}\frac{1}{s^n}Q \log W_n\le \lim_{n\rightarrow \infty} \frac{1}{\theta s^n}\log Q W_n^{\theta}=\lim_{n\to\infty}f_n.
\end{align*}

As a consequence very strong disorder holds if $f_n<0$ for any $f_n$. As a consequence, strong disorder and very strong disorder are equivalent.
\end{proof}

\subsection{Change of measure and environment tilting}

The result of the previous section assures that we can estimate the free energy if we can bound accurately some non integer moment of $W_n$.
Now we present a method to estimate non-integer moment via measure change, it has been introduced to show disorder relevance in the case of wetting on non hierarchical lattice \cite{GLT} and used since in several different contexts since, in particular for directed polymer models on $\Z^d$, \cite{Lac}.
Yet, for the directed polymer on hierarchical lattice, the method is remarkably simple to apply, and it seems to be the ideal context to present it.
\\
Let $\tilde Q$ be any probability measure such that $Q$ and $\tilde Q$ are mutually absolutely continuous.
Using H\"older inequality we observe that 
\begin{equation}
Q W_n^{\theta}
   = \tilde Q \frac{\dd Q}{\dd \tilde Q} W_n^{\theta}
   \le \left[\tilde Q\left(\frac{\dd Q}
     {\dd \tilde Q}\right)^{\frac{1}{1-\theta}}\right]^{(1-\theta)}\left(\tilde Q W_n\right)^{\theta}.
     \label{holder}
\end{equation}
Our aim is to find a measure $\tilde Q$ such that the term $\left[\tilde Q\left(\frac{\dd Q}
     {\dd \tilde Q}\right)^{\frac{1}{1-\theta}}\right]^{(1-\theta)}$ is not very large (i.e.\ of order $1$), and which significantly lowers the expected value of $W_n$. To do so we look for $\tilde Q$ which lowers the value of the environment on each site, by exponential tilting. For $b<s$ it sufficient to lower the value for the environment uniformly of every site of $D_n\setminus \{A,B\}$ to get a satisfactory result, whereas for the $b=s$ case, on has to do an inhomogeneous change of measure. We present the change of measure in a united framework before going to the details with two separate cases.

Recall that $V_i$ denotes the sites of $D_i\setminus D_{i+1}$, and that the number of sites in $D_n$ is

\begin{align}\label{Dnn}
|D_n\setminus \{A,B\}|=\sum_{i=1}^n |V_i|=\sum_{i=1}^n (s-1)b^is^{i-1}=\frac{(s-1)b((sb)^n-1)}{sb-1}
\end{align}

We define $\tilde Q= \tilde Q_{n,s,b}$ to be the measure under which the environment on the site of the $i$-th generation for $i\in\{1,\dots,n\}$ are standard gaussians with mean  $-\delta_i=\delta_{i,n}$, where $\delta_{i,n}$ is to be defined.
The density of $\tilde Q$ with respect to $Q$ is given by

\begin{align*}
\frac{\dd \tilde Q}{\dd Q}(\go)=\exp\left(-\sum^n_{i=1}\sum_{z\in V_i} (\delta_{i,n}\omega_z+\frac{\delta^2_{i,n}}{2})\right).
\end{align*}

\noindent  As each path in $D_n$ intersects $V_i$ on $s^{i-1}(s-1)$ sites, this change of measure lowers the value of the Hamiltonian \eqref{eq: Hn} by $\sum_{i=1}^n s^{i-1}(s-1)\delta_{i,n}$ on any path. Therefore, 
both terms can be easily computed,
%at least in the gaussian case:

\begin{eqnarray}
	\tilde Q\left(\frac{\dd Q}
     {\dd \tilde Q}\right)^{\frac{1}{1-\theta}} =
     \exp \left\lbrace \frac{\theta}{2(1-\theta)} \sum^n_{i=1}|V_i| \delta^2_{i,n} \right\rbrace.
	\label{cost}
\end{eqnarray}

\begin{eqnarray}
	\left(\tilde Q W_n\right)^{\theta} = \exp \left\lbrace -\beta \theta \sum^n_{i=1}  s^{i-1}(s-1) \delta_{i,n}\right\rbrace.
	\label{gain}
\end{eqnarray}

\noindent Replacing \eqref{gain} and \eqref{cost} back into \eqref{holder} gives

\begin{eqnarray}
	u_n \leq \exp \left\lbrace \theta \sum^n_{i=1} \left(  \frac{|V_i|\delta^2_{i,n}}{2(1-\theta)}- \beta s^{i-1}(s-1) \delta_{i,n} \right) \right\rbrace.
	\label{eq:ineq0}
\end{eqnarray}

\vspace{2ex}

\noindent When $\delta_{i,n}= \delta_n$ (i.e. when the change of measure is homogeneous on every site) the last expression becomes simply

\begin{eqnarray}
	u_n \leq \exp \left\lbrace \theta \left(  \frac{|D_n\setminus\{A,B\}|\delta^2_{n}}{2(1-\theta)}- (s^n-1)\beta \delta_{n} \right) \right\rbrace.
	\label{ineq1}
\end{eqnarray}

\noindent In either case, the rest of the proof then consists in finding convenient values for $\delta_{i,n}$ and
$n$ large enough to insure that $(i)$ from Proposition \ref{th:fracmom} holds.

%--------------------------------------------------------------------------------------------------------------------------

%--------------------------------------------------------------------------------------------------------------------------

%--------------------------------------------------------------------------------------------------------------------------

\subsection{Homogeneous shift method: $s>b$}

\begin{proof}[Proof of the left inequality in Theorem \ref{th:bs}, in the gaussian case]

Let $0<\theta<1$ be fixed (say $\theta=1/2$) and $\delta_{i,n}= \delta_n :=(sb)^{-n/2}$.

\noindent Observe from \eqref{Dnn} that $|D_n\setminus\{A,B\}| \delta^2_n \leq 1$, so that \eqref{ineq1} implies

\begin{align*}
u_n\le \exp\left(\frac{\theta}{2(1-\theta)}-\theta \gb (s/b)^{n/2}\frac{s-1}{s}\right).
\end{align*}

\noindent Taking $n=\frac{2(|\log \gb|+\log c_3)}{\log s-\log b}$, we get

\begin{align*}
u_n\le \exp\left(\frac{\theta}{2(1-\theta)}-\frac{\theta c_5s}{s-1}\right).
\end{align*}

\noindent Choosing $\theta=1/2$ and $c_3$ sufficiently large, we have
\begin{align}
f_n=s^{-n}\log a_{\theta}b^{\frac{1-\theta}{s-1}}u_n\le -s^{-n},
\end{align}
so that Proposition \ref{th:fracmom} gives us the conclusion

\begin{align*}
p(\beta)-\lambda(\beta) \le -s^{-n}=-(\gb/c_3)^{\frac{2\log s}{\log s -\log b}}.
\end{align*}
\end{proof}

%--------------------------------------------------------------------------------------------------------------------------

%--------------------------------------------------------------------------------------------------------------------------

%--------------------------------------------------------------------------------------------------------------------------

\subsection{Inhomogeneous shift method: $s=b$} \label{inshm}
One can check that the previous method does not give good enough results for the marginal case $b=s$.
One has to do a change of measure which is a bit more refined and for which the intensity of the tilt in proportional to the Green Function on each site. This idea was used first for the marginal case in pinning model on hierarchical lattice (see \cite{L}).

\begin{proof}[Proof of the left inequality in Theorem \ref{th:ss}, the gaussian case]
This time, we set $\gd_{i,n}:=n^{-1/2}s^{-i}$.
Then (recall \eqref{Dnn}), \eqref{eq:ineq0} becomes

\begin{align*}
u_n\le \exp\left(\frac{\theta}{2(1-\theta)}\frac{s-1}{s}-\theta\gb n^{-1/2}\frac{s-1}{s}\right).
\end{align*}

\noindent Taking $\theta=1/2$ and $n=(c_4/\gb)^2$ for a large enough constant $c_4$, we get that $f_n\le -s^n$  and applying Proposition \ref{th:fracmom}, we obtain

\begin{align*}
p(\beta) - \lambda(\beta) \le -s^{-n}=-s^{-(c_4/\gb)^2}=\exp\left(-\frac{c_4^2\log s}{\gb^2}\right).
\end{align*}
\end{proof}

%--------------------------------------------------------------------------------------------------------------------------

%--------------------------------------------------------------------------------------------------------------------------

%--------------------------------------------------------------------------------------------------------------------------
\subsection{Bounds for the critical temperature}

From Proposition \ref{th:fracmom}, we have that strong disorder holds if
$a_{\theta}< b^{(1-\theta)/(s-1)}$. Taking logarithms, this condition
reads

\begin{eqnarray*}
\lambda(\theta \beta)-\theta \lambda(\beta)< (1-\theta)\frac{ \log b}{s-1}.
\end{eqnarray*}

\noindent We now divide both sides by $1-\theta$ and let $\theta \to 1$. This
proves part $(iii)$ of Proposition \ref{misc}.

\noindent For the case $b>s$, this condition can be improved by the inhomogeneous shifting method; here, we perform it just in the gaussian case. Recall that

\begin{eqnarray}
	u_n \leq \exp \left\lbrace \theta \sum^n_{i=1} \left( \frac{|V_i|\delta^2_{i,n}}{2(1-\theta)}- \beta s^{i-1}(s-1)\delta_{i,n} \right) \right\rbrace.
\end{eqnarray}

\noindent We optimize each summand in this expression taking $\delta_{i,n}=\delta_i = (1-\theta) \beta / b^i$. Recalling
that $|V_i|= (bs)^{i-1}b(s-1)$, this yields

\begin{eqnarray}
\nonumber
	u_n &\leq& \exp \left\lbrace -\theta (1-\theta) \frac{\beta^2}{2} \frac{s-1}{s}  \sum^n_{i=1} \left( \frac{s}{b}\right)^i \right\rbrace\\
	\nonumber
	&\leq& \exp \left\lbrace -\theta (1-\theta) \frac{\beta^2}{2} \frac{s-1}{s}  \frac{s/b-(s/b)^{n+1}}{1-s/b} \right\rbrace.
\end{eqnarray}

\noindent Because $n$ is arbitrary, in order to guaranty strong disorder it is enough to have ( cf. first condition in Proposition \ref{th:fracmom}) for some $\theta\in(0,1)$

\begin{eqnarray}\nonumber
	\theta (1-\theta) \frac{\beta^2}{2} \frac{s-1}{s}  \frac{s/b}{1-s/b}
	> (1-\theta)\frac{\log b}{s-1}+\log a_\theta.
\end{eqnarray}

\noindent In the case of gaussian variables $\log a_\theta=\theta(\theta-1)\gb^2/2$. This is equivalent to

\begin{eqnarray*}
\frac{\beta^2}{2} > \frac{(b-s)\log b}{(b-1)(s-1)}.
\end{eqnarray*}

\noindent This last condition is an improvement of the bound in part $(iii)$ of
Proposition \ref{misc}.

%if

%\begin{eqnarray*}
%\frac{b-s}{s-1}<1,
%\end{eqnarray*}

%\noindent that is, if $b+1<2s$.

%--------------------------------------------------------------------------------------------------------------------------

%--------------------------------------------------------------------------------------------------------------------------

%--------------------------------------------------------------------------------------------------------------------------

\subsection{Adaptation of the proofs for non-gaussian variables}
\begin{proof}[Proof of the left inequality in Theorem \ref{th:bs} and \ref{th:ss}, the general case]
To adapt the preceding proofs to non-gaussian variables, we have to investigate the consequence of exponential tilting on non-gaussian variables. We sketch the proof in the inhomogeneous case $b=s$, we keep $\gd_{i,n}:=s^{-i}n^{-1/2}$.

Consider $\tilde Q$ with density
\begin{align*}
\frac{\dd \tilde Q}{\dd Q}(\go):=\exp\left(-\sum_{i=1}^n\sum_{z\in V_i}\left(\gd_{i,n}\go_z+\gl(-\gd_{i,n})\right)\right),
\end{align*}
(recall that $\gl(x):=\log Q \exp(x\go)$).
The term giving cost of the change of measure is, in this case,
\begin{eqnarray*}
\left[\tilde Q\left(\frac{\dd Q}{\dd \tilde Q}\right)^{\frac{1}{1-\theta}}\right]^{(1-\theta)}&=&\exp\left((1-\theta)\sum_{i=1}^{n}|V_i|\left[\gl\left(\frac{\theta\gd_{i,n}}{1-\theta}\right)+\frac{\theta}{1-\theta}\gl(-\gd_{i,n})\right]\right)\\
&\leq& \exp\left(\frac{\theta}{(1-\theta)}\sum_{i=1}^n|V_i|\gd_{i,n}^2\right)\, \le \, \exp\left(\frac{\theta}{(1-\theta)}\right)
\end{eqnarray*}
Where the inequality is obtained by using the fact the $\gl(x)\sim_0 x^2/2$ (this is a consequence of the fact that $\go$ has unit variance) so that if $\gb$ is small enough, one can bound every $\gl(x)$ in the formula by $x^2$.

We must be careful when we estimate $\tilde Q W_n$. We have
\begin{align*}
\tilde Q W_n=\exp \left(\sum_{i=1}^n (s-1)s^{i-1}\gl (\gb-\gd_{i,n})-\gl(\gb)-\gl(-\gd_{i,n})\right) Q W_n
\end{align*}
By the mean value theorem
\begin{align*}
\gl (\gb-\gd_{i,n})-\gl(\gb)-\gl(-\gd_{i,n})+\gl(0)=-\gd_{i,n}\left(\gl'(\gb-t_0)-\gl'(-t_0)\right)=-\gd_{i,n}\gb\gl''(t_1),
\end{align*}
for some $t_0\in(0,\gd_{i,n})$ and some $t_1\in (\gb,-\gd_{i,n})$. As we know that $\lim_{\gb\to 0}\gl''(\gb)=1$,
when $\gd_i$ and $\gb$ are small enough, the right-hand side is less than $-\gb\gd_{i,n}/2$.
Hence,
\begin{align*}
\tilde Q W_n\le \exp \left(\sum_{i=1}^n(s-1)s^{i-1}\frac{\gb\gd_{i,n}}{2}\right).
\end{align*}
We get the same inequalities that in the case of gaussian environment, with different constants, which do not affect the proof. The case $b<s$ is similar.
\end{proof}

\section{Fluctuation and localisation results}\label{flucloc}

In this section we use the shift method we have developed earlier to prove fluctuation results

\subsection{Proof of Proposition \ref{fluctuations}}

The statement on the variance is only a consequence of the second statement. Recall that the random variable $\go_z$ here are i.i.d. centered standard gaussians, and that the product law is denoted by $Q$.
We have to prove
\begin{equation}
Q\left\{ \log Z_n \in [a, a+\gb\gep(s/b)^{n/2}]\right\}\le 4\gep \quad \forall \gep>0, n\ge 0, a\in \bbR \label{fluc}
\end{equation}
Assume there exist real numbers $a$ and $\gep$, and an integer $n$ such that \eqref{fluc} does not hold, i.e.

\begin{equation}
Q\left\{ \log \bar Z_n \in [a, a+\gb\gep(s/b)^{n/2})\right\}> 4\gep.
\end{equation}
Then one of the following holds

\begin{equation}\begin{split}\label{abcd}
Q\left\{ \log Z_n \in [a, a+\gb\gep(s/b)^{n/2})\right\}\cap\left\{\sum_{z\in D_n} \go_z\ge 0\right\}&> 2\gep,\\
Q\left\{ \log  Z_n \in [a, a+\gb\gep(s/b)^{n/2})\right\}\cap\left\{\sum_{z\in D_n} \go_z\le 0\right\}&> 2\gep.
\end{split}\end{equation}
We assume that the first line is true. We consider the events related to $Q$ as sets of environments $(\go_z)_{z\in D_n\setminus \{A,B\}}$. We define

\begin{equation}
A_\gep=\left\{ \log Z_n \in [a, a+\gb\gep b^{-n/2})\right\}\cap\left\{\sum_{z\in D_n} \go_z\ge 0\right\},
\end{equation}
and

\begin{equation}
A^{(i)}_\gep=Q\left\{ \log Z_n \in [a-i\gb\gep(s/b)^{n/2}, a-(i-1)\gb\gep(s/b)^{n/2})\right\}.
\end{equation}
Define $\gd=\frac{s^{n/2}}{(s^n-1)b^{n/2}}$.  We define the measure $\tilde Q_{i,\gep}$ with its density:
\begin{equation}
\frac{\dd \tilde Q_{i,\gep}}{\dd Q}(\go):=\exp\left(\left[i\gep\gd^2 \sum_{z\in D_n} \go_z\right]-\frac{i^2\gep^2\gd^2|D_n\setminus\{A,B\}|}{2}\right).
\end{equation}
If the environment $(\go_z)_{z\in D_n}$ has law $Q$ then $(\hat \go_z^{(i)})_{z\in D_n}$ defined by
\begin{equation}
 \hat \go_z^{(i)}:= \go_z+\gep i \gd,
\end{equation}
 has law $ \tilde Q_{i,\gep}$. Going from $\go$ to $\hat \go^{(i)}$, one increases the value of the Hamiltonian by $\gep i(s/b)^{n/2}$ (each path cross $s^n-1$ sites).
Therefore if $(\hat\go^{(i)}_z)_{z\in D_n}\in A_{\gep}$, then $(\go_z)_{z\in D_n}\in A^{(i)}_{\gep}$.
From this we have $\tilde Q_{i,\gep} A_\gep \le Q A^{(i)}_{\gep}$, and therefore
\begin{equation}
Q A^{(i)}_{\gep}\, \ge \int_{A_\gep} \frac{\dd \tilde Q_{i,\gep}}{\dd Q}Q(\dd \go)\ge \exp(-(\gep i)^2/2)Q(A_\gep).
\end{equation}
The last inequality is due to the fact that the density is always larger than $\exp(-(\gep i)^2/2)$ on the set $A_\gep$ (recall its definition and the fact that $|D_n\setminus\{A,B\}|\gd^2\le 1$).
Therefore, in our setup, we have 
\begin{equation}
Q A^{(i)}_{\gep}> \gep, \quad  \forall i\in[0,\gep^{-1}].
\end{equation}
As the $A^{(i)}_{\gep}$ are disjoints, this is impossible.
If we are in the second case of \eqref{abcd}, we get the same result by shifting the variables in the other direction. \qed

\subsection{Proof of Proposition \ref{fluc2}}

Let us suppose that there exist $n$, $\gep$ and $a$ such that
\begin{equation}
Q\left\{ \log Z_n \in [a, a+\gb\gep\sqrt{n})\right\}> 8\gep.
\end{equation}
We define $\delta_{i,n}=\delta_i:=\gep s^{1-i}(s-1)^{-1}n^{-1/2}$.
Then one of the following inequality holds (recall the definition of $V_i$)

\begin{equation}\begin{split} \label{fghi}
Q\left\{ \log Z_n \in [a, a+\gb\gep\sqrt{n})\right\}\cap\left\{\sum_{i=1}^{n}\gd_i\sum_{z\in V_i}\go_z\ge 0\right\}&> 4\gep,\\
Q\left\{ \log Z_n \in [a, a+\gb\gep\sqrt{n})\right\}\cap\left\{\sum_{i=1}^{n}\gd_i\sum_{z\in V_i}\go_z\le 0\right\}&> 4\gep.
\end{split}\end{equation}
We assume that the first line holds and define

\begin{equation}
A_{\gep}=\left\{ \log Z_n \in [a, a+\gb\gep\sqrt{n})\right\}\cap\left\{\sum_{i=1}^{n}\gd_i\sum_{z\in V_i}\go_z\ge 0\right\}
\end{equation}
And
\begin{equation}
A_{\gep}^{(j)}=\left\{ \log Z_n \in [a-j\gb\gep\sqrt{n}, a-(j-1)\gb\gep\sqrt{n})\right\}
\end{equation}

\begin{equation}
j\gb\sum_{i=1}^n\gd_i (s-1)s^{i-1}=j\gb\gep \sqrt{n}.
\end{equation}
Therefore, an environment $\go\in A_{\gep}$ will be transformed in an environment in $A^{(j)}_{\gep}$.

We define $\tilde Q_{j,\gep}$ the measure whose Radon-Nicodyn derivative with respect to $Q$ is
\begin{equation}
\frac{\dd \tilde Q_{i,\gep}}{\dd Q}(\go):=\exp\left(\left[j\sum_{i=1}^n\gd_i \sum_{z\in V_i} \go_z\right]-\sum_{i=1}^{n}\frac{j^2\gd_i^2|V_i|}{2}\right).
\end{equation}
We can bound the deterministic term.
\begin{equation}
\sum_{i=1}^{n}\frac{j^2\gd_i^2|V_i|}{2}=j^2\gep^2 \sum_{i=1}^n \frac{s}{2(s-1)}\le j^2\gep^2.
\end{equation}
For an environment $(\go_z)_{z\in D_n\setminus\{A,B\}}$, define $(\hat \go_z^{(j)})_{z\in D_n\setminus\{A,B\}}$ by
\begin{equation}
 \hat \go_z^{(j)}:= \go_z+j\gep\gd_i, \quad   \forall z \in V_i.
\end{equation}
If $(\go_z)_{z\in D_n\setminus\{A,B\}}$ has $Q$, then $(\hat \go_z^{(j)})_{z\in D_n\setminus\{A,B\}}$ has law $\tilde Q_{j,\gep}$.
When one goes from $\go$ to $\hat\go^{(j)}$, the value of the Hamiltonian is increased by
\begin{equation*}
 \sum_{i=1}^n j\gd_i s^{i-1}(s-1)=\gep \sqrt{n}.
\end{equation*}
Therefore, if $\hat \go^{(j)}\in A_{\gep}$, then $\go\in A^{(j)}_{\gep}$, so that 
\begin{equation*}
 Q A^{(j)}_{\gep}\ge \tilde Q_{j,\gep} A_{\gep}.
\end{equation*}
Because of the preceding remarks
\begin{equation}
Q A^{(j)}_{\gep}\ge \tilde Q_{j,\gep} A_{\gep}=\int_{A_{\gep}}\frac{\dd \tilde Q_{i,\gep}}{\dd Q}Q(\dd\go) \ge \exp\left(-j^2\gep^2\right)Q A_{\gep}.
\end{equation}
The last inequality comes from the definition of $A_{\gep}$ which gives an easy lower bound on the Radon-Nicodyn derivative.
For $j\in [0,(\gep/2)^{-1}]$, this implies that $Q A^{(j)}_{\gep}> 2\gep$. As they are disjoint events this is impossible. The second case of \eqref{fghi} can be dealt analogously. \qed

\subsection{Proof of Corollary \ref{locloc}}

Let $g\in \gG_n$ be a fixed path. For $m\ge n$, define
\begin{equation}
 Z_m^{(g)}:=\sum_{\{ g'\in \gG_m: g|_n=g\}}\exp\left(\gb H_m(g')\right).
\end{equation}
With this definition we have
\begin{equation}
 \mu_m(\gamma|_n=g)=\frac{Z_m^{(g)}}{Z_m}.
\end{equation}
To show our result, it is sufficient to show that for any constant $K$ and any distinct $g,g'\in \gG_n$
\begin{equation}
 \lim_{m\to\infty}Q\left( \frac{\mu_m(\gamma|_n=g)}{\mu_m(\gamma|_n=g')}\in[K^{-1},K]\right)=0.
\end{equation}
For $g$ and $g'$ distinct, it is not hard to see that
\begin{equation}
 \log\left( \frac{\mu_m(\gamma|_n=g)}{\mu_m(\gamma|_n=g')}\right)=\log Z_m^{(g)}-\log Z_m^{(g')}=: \log Z^{(0)}_{m-n}+X,
\end{equation}
where $Z^{(0)}_{m-n}$ is a random variable whose distribution is the same as the one of $Z_{m-n}$, and $X$ is independent of $Z^{(0)}_{m-n}$.
We have
\begin{multline}
 Q\left(\log \left(\frac{\mu_m(\gamma|_n=g)}{\mu_m(\gamma|_n=g')}\right)\in [-\log K, \log K]\right)\\
=Q \left[ Q\left(\log Z^{(0)}_{m-n}\in [-\log K -X, \log K -X] \, \big| \, X\right)\right]
 \\
\le \max_{a\in \R} Q \left(\log Z_{m-n}\in [a,a+2\log K]\right).
\end{multline}
Proposition \ref{fluctuations} and \ref{fluc2} show that the right--hand side tends to zero. \qed

%-------------------------------------------------------------------------------------------------------------------------------

%---------------------------------------------------------------------WEAK DISORDER MEASURE-------------------------------------

%-------------------------------------------------------------------------------------------------------------------------------

%-------------------------------------------------------------------------------------------------------------------------------

\section{The weak disorder polymer measure}\label{weakpol}

Comets and Yoshida introduced in \cite{CY} an infinite
volume Markov chain at weak disorder that corresponds in some sense to the limit of the polymers
measures $\mu_n$ when $n$ goes to infinity. We perform the same construction here.
The notation is more cumbersome in our setting.

Recall that $\Gamma_n$ is the space of directed paths from $A$ to $B$ in $D_n$. Denote by $P_n$ the uniform law on $\gG_n$. For $g\in \Gamma_n$,
$0\leq t \leq s^n-1$, define $W_{\infty}(g_t,g_{t+1})$ by performing the same construction that leads to $W_{\infty}$, but taking $g_t$ and $g_{t+1}$ instead of $A$ and $B$ respectively. On the classical directed polymers on $\Z^d$, this would be equivalent to take the $(t,g_t)$ as the initial point of the polymer.

We can now define the weak disorder polymer measure for $\beta < \beta_0$. We define $\gG$ as the projective limit of $\gG_n$ (with its natural topology), the set of path on $D:=\bigcup_{n\ge 1} D_n$. As for finite path, we can define, for $\bar g\in \gG$, its projection onto $\gG_n$, $\bar g|_n$.
We define

\begin{eqnarray}\label{wdpol}
\mu_{\infty}(\bar \gamma|_n =g):= \frac{1}{W_{\infty}} \exp \{ \beta H_n(g) - (s^n-1) \lambda(\beta) \} \prod^{s^n-1}_{i=0} W_{\infty}({g}_i,{g}_{i+1}) \, P_n({\bar \gamma|_n = g}).
\end{eqnarray}

Let us stress the following:

\begin{itemize}
      
\item Note that the projection on the different $\Gamma_n$ are consistent (so that our definition makes sense)
  
      \begin{eqnarray*}
      \mu_{\infty}(\bar \gamma|_n = g)= \mu_{\infty}\left((\bar \gamma|_{n+1})|_n = g\right).
      \end{eqnarray*}

\item Thanks to the martingale convergence for both the numerator and the
      denominator, for any ${\bf s} \in \Gamma_n$,
      
      \begin{eqnarray*}
      \lim_{k\to +\infty} \mu_{k+n}(\gamma|_n = g) = \mu_{\infty}(\bar \gamma|_n = g).
      \end{eqnarray*}
Therefore, $\mu_{\infty}$ is the only reasonable definition for the limit of $\mu_n$.
%\item This construction allows us to define a convenient sequence of (random) paths $\gamma^{(n)}\in \Gamma_n$,
%      such that $\gamma^{(n+1)}|_n=\gamma^{(n)}$: chose $\gamma^{(1)}$ at random from $\Gamma_1$
%      following $\mu_{\infty}$. Then, for $i=0,...,s-1$ choose a path linking $\gamma^{(1)}(i)$ to $\gamma^{(1)}(i+1)$
%      according to the following law: for any path ${\bf s}$ of length $s$ linking $\gamma^{(1)}(i)$ to $\gamma^{(1)}(i+1)$,
      
%      \begin{multline}
%      \mu_{\infty,\gamma_i,\gamma_{i+1}}(\gamma'={\bf s})
%      =\frac{1}{W_{\infty}(\gamma_i,\gamma_{i+1})}\exp{\beta H_1({\bf s})-(s-1)\lambda(\beta)}
%       \prod^{s-1}_{j=0} W_{\infty}({\bf s}_j,{\bf s}_{j+1}) \, P_1(\gamma' = {\bf s}),
%      \end{multline}

%      \noindent where it is understood that $H_1$ is defined using the random variables of generation $2$ that lies between $\gamma^{(1)}(i)$ and $\gamma^{(1)}(i+1)$.

%      \noindent This procedure is easily iterated to obtain the wanted sequence of paths.
%      It is easy to see that the law of $\gamma^{(n)}$ is the law of a path chosen randomly
%      from $\Gamma_n$ according to $\mu_{\infty}$.
\end{itemize}

\vspace{3ex}

It is an easy task to prove the law of large numbers for the time-averaged quenches mean of the energy.
This follows as a simple consequence of the convexity of $p(\beta)$.

\begin{proposition}
At each point where $p$ admits a derivative,
 \begin{eqnarray*}
   \lim_{n\to +\infty} \frac{1}{s^n} \mu_{n}(H_n(\gamma)) \to  p'(\beta), \quad Q-a.s..
 \end{eqnarray*}
\end{proposition}
\begin{proof}
 It is enough to observe that 
 \begin{eqnarray*}
  \frac{d}{d \beta} \log Z_n = \mu_n(H_n(\gamma)),
 \end{eqnarray*}
 
 \noindent then use the convexity to pass to the limit.
\end{proof}

\vspace{2ex}

We can also prove a quenched law of large numbers under our infinite volume
measure $\mu_{\infty}$, for almost every environment. The proof is very easy,
as it involves just a second moment computation. 

\begin{proposition}
 At weak disorder,
 \begin{eqnarray*}
   \lim_{n\to +\infty}\frac{1}{s^n}H_n(\bar \gamma|_n)=\lambda'(\beta),\quad \mu_{\infty}-a.s., \, Q-a.s..
 \end{eqnarray*}
\end{proposition}

\begin{proof}
We will consider the following auxiliary measure (size biased measure) on the environment

\begin{eqnarray*}
\overline{Q}(f(\go)) = Q(f(\go) W_{+\infty}).
\end{eqnarray*}

\noindent So, $Q$-a.s. convergence will follow from $\overline{Q}$-a.s. convergence.
This will be done by a direct computation of second moments. Let us write 
$\Delta = Q(\omega^2 e^{\beta \omega - \lambda(\beta)})$.

\begin{eqnarray*}
&\ &\overline{Q} \left( \mu_{\infty}(|H_n(\bar \gamma|_n)|^2) \right)\\
&=& Q \left[ P_n\left( |H_n(\gamma)|^2 \exp \{ \beta H_n(\gamma)- (s^n-1) \lambda(\beta)\} \prod^{s^n-1}_{i=0} W_{\infty}(\gamma_i,\gamma_{i+1})\right)\right]\\
&=& Q \left[ P_n( |H_n(\gamma)|^2 \exp \{ \beta H_n(\gamma)- (s^n-1) \lambda(\beta)\})\right]\\
&=& Q \left[ P_n( |\sum^{s^n}_{t=1}\omega(\gamma_t)|^2 \exp \{ \beta H_n(\gamma)- (s^n-1) \lambda(\beta)\})\right]\\
&=& Q\left[ \sum^{s^n-1}_{t=1}P_n\left(|\omega(\gamma_t)|^2\exp \{ \beta H_n(\gamma)- (s^n-1) \lambda(\beta)\}\right)\right]\\
&\ &+ Q \left[\sum_{1\le t_1\neq t_2\le s^{n}-1}P_n\left(\omega(\gamma_{t_1})\omega(\gamma_{t_2})\exp \{ \beta H_n(\gamma)- (s^n-1) \lambda(\beta)\}\right)\right]\\
&=& (s^n-1) \Delta+(s^n-1)(s^n-2)(\lambda'(\beta))^2,
\end{eqnarray*}

\noindent where we used independence to pass from line two to line three. So, recalling that $\overline{Q} ( \mu_{\infty}(H_n(\bar \gamma_n))= (s^n-1) \lambda'(\beta)$, we
have

\begin{eqnarray*}
&\ &\overline{Q} \left( \mu_{\infty}(|H_n(\gamma^{(n)}) - (s^n-1) \lambda'(\beta)|^2) \right)\\
&=& (s^n-1) \Delta+(s^n-1)(s^n-2)(\lambda'(\beta))^2 
- 2 (s^n-1) \lambda'(\beta)\overline{Q} \left( \mu_{\infty}(H_n(\bar \gamma_n)) \right)\\
&\ & \quad + \, (s^{n}-1)^2(\lambda'(\beta))^2\\
%&=& s^n \Delta+s^n(s^n-1)(\lambda'(\beta))^2 
%- 2 s^n \lambda'(\beta) s^n \lambda'(\beta)
%+ s^{2n}(\lambda'(\beta))^2\\
%&=& s^n \Delta+s^n(s^n-1)(\lambda'(\beta))^2 
%- 2 s^{2n} (\lambda'(\beta))^2 
%+ s^{2n}(\lambda'(\beta))^2\\
&=& (s^n-1) \left( \Delta - (\lambda'(\beta)^2)\right).
\end{eqnarray*}

\noindent Then

\begin{eqnarray*}
\overline{Q}  \mu_{\infty}\left( \left|\frac{H_n(\bar\gamma_n) - (s^n-1) \lambda'(\beta)}{s^n} \right|^2\right) 
\leq \frac{1}{s^n}\left( \Delta - (\lambda'(\beta)^2)\right),
\end{eqnarray*}

\noindent so the result follows by Borel-Cantelli.
\end{proof}

\section{Some remarks on the bond--disorder model}\label{bddis}

In this section, we shortly discuss, without going through the details, how the methods we used in this paper could be used (or could not be used) for the model of directed polymer on the same lattice with disorder located on the bonds.
\\

In this model to each bond $e$ of $D_n$ we associate i.i.d.\ random variables $\go_e$.
We consider each set $g\in \gG_n$ as a set of bonds and define the Hamiltonian as
\begin{equation}
H^{\omega}_n(g) = \sum_{e\in g} \go_e,
\end{equation}

\noindent  The partition function $Z_n$ is defined as
\begin{equation}
 Z_n:= \sum_{g\in \gG_n} \exp(\gb H_n(g)).
\end{equation}
One can check that is satisfies the following recursion 

\begin{equation}\begin{split}
Z_0& \, \stackrel{ \mathcal L}{=}\, \exp(\gb\go)\\
Z_{n+1}&\, \stackrel{\mathcal L)}{=}\, \sum_{i=1}^{b}Z_n^{(i,1)}Z_n^{(i,2)}\dots Z_n^{(i,s)}.
\end{split}
\end{equation}

\noindent where equalities hold in distribution and and $Z_n^{i,j}$ are i.i.d.\ distributed copies of $Z_n$.
Because of the loss of the martingale structure and the homogeneity of the Green function in this model
(which is equal to $b^{-n}$ on each edge), Lemma \ref{martingale} does not hold, and we cannot prove part $(iv)$ in Proposition \ref{misc}, Theorem \ref{th:ss} and Proposition \ref{fluc2} for this model. Moreover we have to change $b\le s$ by $b<s$ in $(v)$ of Proposition \ref{misc}.
Moreover, the method of the control of the variance would give us a result similar to \ref{th:ss} in this case
\begin{proposition}\label{ssbd}
When $b$ is equal to $s$, on can find constants $c$ and $\gb_0$ such that for all $\gb\le \gb_0$
\begin{equation}
0\le \lambda(\beta)-p(\beta) \le \exp\left(-\frac{c}{\gb^2}\right).
\end{equation}
\end{proposition}
However, we would not be able to prove that annealed and free energy differs at high temperature for $s=b$ using our method. The techniques used in \cite{GLT_marg} or \cite{Lac} for dimension $2$ should be able to tackle this problem, and show marginal disorder relevance in this case as well. 
\medskip

\noindent {\bf Acknowledgements:} The authors are very grateful to Francis Comets and Giambattista Giacomin for their suggestion to work on this subject and many enlightning discussions. This work was partially supported by CNRS, UMR $7599$ ``Probabilit\'es et
Mod\`eles Al\'eatoires". H.L.\ acknowledges the support of ANR grant POLINTBIO. G.M.\ acknowledges the support of Beca Conicyt-Ambassade de France.


\begin{thebibliography}{99}

\bibitem{BZ} P. M. Bleher and E. Zhalis, {\it Limit Gibbs distributions for the Ising model on hierarchical lattices}, (Russian)  Litovsk. Mat. Sb. {\bf 28} 2  (1988), 252-268;  translation in  Lithuanian Math. J.  {\bf 28} 2  (1988),  127--139 

\bibitem{Birk} M. Birkner, {\it A Condition for Weak Disorder for Directed Polymers in Random Environment},  Elec. Comm. Probab {\bf 9} (2004), 22--25.

\bibitem{BS} M. Birkner and R. Sun, {\it Annealed vs Quenched Critical Points for a Random Walk Pinning Model}, to appear in Ann. Inst. H. Poinc., Probab. and Stat., arXiv:0807.2752v2 [math.PR].

\bibitem{B} E. Bolthausen, {\it A note on diffusion of directed polymeres in a random environment}, Commun. Math. Phys.  {\bf 123} (1989), 529-534.

\bibitem{Buf} E. Buffet, A. Patrick and J.V. Pule,  {\it Directed polymers on trees: a martingale approach},  Journal Of Physics {\bf 26} (1993), 1823-1834.

\bibitem{CC} A. Camanes and P. Carmona, {\it The critical temperature of a Directed Polymer in a random environment}, to appear in Markov Proc. Relat. Fields.

\bibitem{CH_al} P. Carmona and Y. Hu, {\sl Strong disorder implies strong localization for directed
polymers in a random environment}, ALEA {\bf 2} (2006), 217--229.

\bibitem{CH_ptrf} P. Carmona and Y. Hu,
{\sl On the partition function of a directed polymer in a random Gaussian environment} ,
Probab. Theor. Relat. Fields {\bf 124} 3 (2002), 431-457.

\bibitem{Cnotes} Comets, F. {\sl Unpublished lecture notes}, (2008).

\bibitem{CPV} F. Comets, S. Popov and M. Vachkovskaia, {\it The number of open paths in an oriented $\rho$-percolation model},  J. Stat. Phys {\bf 131} (2008),  357--379.

\bibitem{CSY} F. Comets ,T. Shiga and N. Yoshida, {\it Directed Polymers in a
   random environment: strong disorder and path localization}, Bernouilli {\bf 9} 4 (2003), 705-723 .

\bibitem{CSY_rev} F. Comets, T. Shiga, and N. Yoshida {\sl Probabilistic Analysis of Directed Polymers in a Random Environment: a Review }, Adv. Stud. Pure Math. {\bf 39} (2004), 115--142 .

\bibitem{CV} F. Comets and V. Vargas {\sl Majorizing multiplicative cascades for directed polymers in random media},  ALEA  {\bf 2} (2006), 267--277 

\bibitem{CY_cmp}
F. Comets and  N. Yoshida {\sl Brownian directed polymers in random environment},  Commun. Math. Phys. {\bf 254}  2 (2005), 257--287.

\bibitem{CY} F. Comets and N. Yoshida {\it Directed polymers in a random environment
   are diffusive at weak disorder}, Ann. Probab. {\bf 34} 5 (2006), 1746--1770 .

\bibitem{CD} J. Cook and  B. Derrida  {\it Polymers on Disordered Hierarchical Lattices: A
   Nonlinear Combination of Random Variables}, J. Stat. Phys. {\bf 57} 1/2 (1989), 89--139.

\bibitem{DF} B. Derrida and H. Flyvbjerg, {\it A new real space renormalization and its Julia set}, J. Phys. A:
Math. Gen. {\bf 18} (1985), L313-L318

\bibitem{DG} B. Derrida and E. Gardner, {\sl Renormalisation group study of a disordered model},
J. Phys. A: Math. Gen. \textbf{17} (1984), 3223--3236.

\bibitem{DGr} B. Derrida and R.B. Griffith,  {\it Directed polymers on disordered hierarchical lattices},
Europhys. Lett.  {\bf 8} 2 (1989), 111--116.

\bibitem{DHV} B. Derrida, V. Hakim and J. Vannimenius,
  \textit{Effect of disorder on two-dimensional wetting}, J. Stat.
  Phys. {\bf 66} (1992), 1189--1213.

\bibitem{DGLT}  B. Derrida, G. Giacomin, H. Lacoin and F.L. Toninelli, {\sl Fractional moment bounds and disorder relevance for pinning models},  Commun. Math. Phys. {\bf 287} (2009), 867--887.

\bibitem{Fr} J Franchi, {\it Chaos multiplicatif : un traitement simple et complet de la fonction de partition}, S\'eminaire de Probabilit\'es de Strasbourg {\bf 29} (1995), 194--201. 

\bibitem{GLT} G. Giacomin, H. Lacoin and F. L. Toninelli, \emph{
    Hierarchical pinning models, quadratic maps and quenched
    disorder},  to appear in Probab. Theor. Rel. Fields, arXiv:0711.4649 [math.PR].


\bibitem{GLT_marg}  G. Giacomin, H. Lacoin and F.L. Toninelli, {\sl Marginal relevance of disorder for pinning models}, to appear in Comm. Pure Appl. Math., arXiv:0811.4723v1 [math.ph].

\bibitem{GM_h} T. Garel and C. Monthus,  {\sl Critical points of
    quadratic renormalizations of random variables and phase
    transitions of disordered polymer models on diamond lattices},   Phys. Rev. E {\bf 77} (2008) 021132. 

\bibitem{Ham} B.M. Hambly and J.H. Jordan {\it A random hierarchical lattice: the series-parallel graph and its properties}, Adv. Appl. Prob., {\bf 36} (2004) 824-838. 


\bibitem{HK} B.M. Hambly and T Kumagai, {\sl Diffusion on the scaling limit of the critical percolation
      cluster in the diamond hierarchical lattice}, (2008) preprint .

\bibitem{HH} C.L. Henley and D.A. Huse, {\sl Pinning and roughening of domain wall in Ising systems due to random
       impurities}, Phys. Rev. Lett.  {\bf 54} (1985) 2708--2711.

\bibitem{IS} J.Z. Imbrie and T. Spencer, {\it Diffusion of Directed Polymers in a Random Environment},  Jour. Stat. Phys. {\bf 52} (1988), 609--622. 

\bibitem{J} K. Johansson, {\it Transversal fluctuation for increasing subsequences on the plane}, Probab. Theor. Rel. Fields {\bf 116} (2000), 445--456.

\bibitem{Kad} L.P. Kadanoff, {\it Notes on Migdals recursion formulae}, Ann. Phys. {\bf 100} (1976), 359--394.

\bibitem{KP} J.P. Kahane and J. Peyri\`ere, {\it Sur certaines martingales de Benoit Mandelbrot}, Adv. Math. {\bf 22} (1976), 131--145.

\bibitem{KG1} M. Kaufmann and R.B. Griffith, {\it Spin systems on hierarchical lattices. Introduction and thermodynamic limit},  Phys. Rev. B {\bf 3}  26  (1982), no. 9, 5022--5032.

\bibitem{KG2} M. Kaufmann and R.B. Griffith, {\it Spin systems on hierarchical lattices. II. Some examples of soluble models},  Phys. Rev. B {\bf 3}  30  (1984),  no. 1, 244--249

\bibitem{L} H. Lacoin, {\sl Hierarchical pinning model with site disorder: Disorder is marginally relevant}, to appear in Probab. Theor. Rel. Fields, arXiv:0807.4864 [math.PR].

\bibitem{Lac} H. Lacoin {\it New bounds for the free energy of directed polymers in dimension $1+1$ and $1+2$}, preprint (2009)  arXiv:0901.0699v1 [math.ph].

\bibitem{LV} Lesigne, Voln\'y  {\it Large deviations for martingales}, Stoch. Proc. Appl.
   {\bf 96}  (2001), 143--159.

\bibitem{Migd} A.A. Migdal, {\it Recurrence equations in gauge field theory}, JETF {\bf 69} (1975), 810--822, 1457--1467.


\bibitem{R_al} S. Roux, A. Hansen, L.R. Da Silva, L.S. Lucena and R.B. Pandey, {\it Minimal path on the hierarchical diamond lattice}, J. Stat. Phys. {\bf 65} 1/2 (1991), 183--204. 

\bibitem{W} J. Wehr, {\it A strong law of large numbers for iterated functions of independent random variables}, J. Stat. Phys {\bf 86} 5/6 (1997), 1373--1384.




\end{thebibliography}
 \end{document}